\documentclass[11pt, letterpaper]{amsart}
\usepackage[utf8]{inputenc}
\usepackage{setspace}
\usepackage{mathrsfs}
\usepackage{amssymb}
\usepackage{latexsym}
\usepackage{amsfonts}
\usepackage{amsmath}
\usepackage{eucal}
\usepackage{bm}
\usepackage{bbm}
\usepackage{graphicx}
\usepackage[english]{varioref}
\usepackage[nice]{nicefrac}
\usepackage[all]{xy}
\usepackage{amsthm}
\usepackage{amssymb,amsthm,upref,amscd}

\def\m{\mathbb{M}}
\def\op{\operatorname}

\def\mmod{\kern-1pt\operatorname{-mod}}

\def\bg{{\bf G}}
\def\bb{{\bf B}}
\def\bt{{\bf T}}
\def\bu{{\bf U}}
\def\bk{\Bbbk}
\def\one{{\bf 1}}
\def\la{\lambda}
\def\bc{{\bf C}}
\def\f{\mathbb{F}}
\def\sn{s{\bf 1}_\theta}

\newtheorem{theorem}{Theorem}[section]
\newtheorem{lemma}[theorem]{Lemma}

\newtheorem{proposition}[theorem]{Proposition}

\theoremstyle{proposition}

\numberwithin{equation}{section}

\begin{document}

\title[extensions of modules for reductive groups]{On the extensions of certain representations of reductive algebraic groups with Frobenius maps}

\author{Xiaoyu Chen}
\address{Department of Mathematics, Shanghai Normal University,
100 Guilin Road, Shanghai 200234, PR China.}
\email{chenxiaoyu@shnu.edu.cn}

\author{Junbin Dong}
\address{Institute of Mathematical Sciences, ShanghaiTech University, 393 Middle Huaxia Road, Pudong, Shanghai 201210, PR China.}
\email{dongjunbin@shanghaitech.edu.cn}


\date{}


\begin{abstract}
Let ${\bf G}$ be a connected reductive algebraic group defined over the finite field $\mathbb{F}_q$ with $q$ elements,  where $q$ is a power of a prime number $p$. Let $\Bbbk$ be a field and we study the extensions of certain  $\bk\bg$-modules in this paper. We show that the extensions of any modules in $\mathscr{O}(\bg)$ by a finite-dimensional $\bk\bg$-module is zero if $p\ne \op{char}\bk\ge5$ or $\op{char}\bk=0$, where $\mathscr{O}(\bg)$ is the principal representation category defined in \cite{D1}. We determine the necessary and sufficient condition for the vanishing of extensions between naive  induced modules. As an application, we give the condition of the vanishing of extensions between simple modules in $\mathscr{O}({\bf G})$ for $\bg=SL_2(\bar{\mathbb{F}}_q)$.
\end{abstract}
\maketitle

\section{Introduction}
Let ${\bf G}$ be a connected reductive algebraic group defined over the finite field $\mathbb{F}_q$ with $q$ elements, where $q$ is a power of a prime number $p$.
Let $\Bbbk$ be a field. All representations of ${\bf G}$ we consider in this paper are over $\Bbbk$. Let $\bf T$ be a maximal  torus contained in a Borel subgroup $\bf B$ and $\theta$  be a character of ${\bf T}$. Thus $\theta$  can also be regarded as a character of ${\bf B}$ by letting ${\bf U}$ (the unipotent radical of ${\bf B}$) act trivially. Motivated by the study of the Verma module in the representations of complex semisimple Lie algebras, N. H. Xi defined the naive induced module $\mathbb{M}(\theta)=\Bbbk{\bf G}\otimes_{\Bbbk{\bf B}}{\Bbbk}_\theta$ in \cite{Xi} and it was deeply studied in \cite{CD1}.  This induced module has a composition series of finite length  if $\op{char}\Bbbk\neq \op{char} \bar{\mathbb{F}}_q$. Thus we get a large class of simple $\Bbbk {\bf G}$-modules. When $\Bbbk =\bar{\mathbb{F}}_q$, the situation is much more complicated.  In this case, the  induced module $\mathbb{M}(\theta)$ may have infinite length.  The paper \cite{C1} gives a necessary and sufficient condition for the induced module $\mathbb{M}(\theta)$ having a composition series  of finite length. Moreover, X.Y. Chen showed that $\mathbb{M}(\theta)$ always has a simple head (see \cite{C2}), and classified and constructed all the abstract simple modules with $ {\bf B}$-stable lines.

In the paper \cite{D1}, J.B. Dong introduced the principal representation category $\mathscr{O}({\bf G})$, which  was conjectured to be a highest weight category in the sense of \cite{CPS} when $\Bbbk =\mathbb{C}$.  In particular, the extensions of the simple modules in $\mathscr{O}({\bf G})$ may have good properties.  But soon after,  X.Y. Chen constructed a counter example (see \cite{C3}) to show that this conjecture is not true.  His example also tells us  that the  category $\mathscr{O}({\bf G})$ is more complicated than we thought even for ${\bf G}= SL_2(\overline{\mathbb{F}}_q)$. So we are interested in the extensions of the modules in the category $\mathscr{O}({\bf G})$.  In this paper, we give the necessary and sufficient condition for  the extensions  between certain modules in  $\mathscr{O}({\bf G})$ being zero. However, when  $\op{Ext}_{\bk\bg}^1(M,N)\ne 0$ for $M,N \in \mathscr{O}({\bf G})$,  the extensions may be very complicated. We prove that $\op{Ext}_{\bk\bg}^1(M,N)\ne 0$ by constructing a specific non-split extension between $M$ and $N$.

 This paper is organized as follows:  In Section 2, we give some notations and preliminary results. In Section 3, we show that $\op{Ext}_{\bk\bg}^1(M,N)=0$ for $M, N\in\mathscr{O}(G)$ with $\dim N<\infty$ if $p\ne \op{char}\bk\ge5$ or $\op{char}\bk=0$. In Section 4, we show that for any field $\bk$ and two characters $\lambda,\mu$ of $\bf T$,  $\op{Ext}_{\bk\bg}^1(\m(\lambda),\m(\mu))=0$ if and only if $\lambda|_{\bf C}\ne\mu|_{\bf C}$, where $\bc$ is the center of $\bg$. In Section 5, we determine the condition  for the vanishing of  extensions of simple modules in $\mathscr{O}({\bf G})$ for ${\bf G}= SL_2(\bar{\mathbb{F}}_q)$ when
$p \neq \op{char}\Bbbk\geq 3$ or $\op{char}\bk=0$.

\section{Preliminaries}

 As in the introduction,  ${\bf G}$ is a connected reductive algebraic group defined over $\mathbb{F}_q$ with the standard Frobenius homomorphism $\text{Fr}$ induced by the automorphism $x\mapsto x^q$ on $\bar{\mathbb{F}}_q$, where $q$ is a power of a prime number $p$. Let ${\bf B}$ be an $\text{Fr}$ -stable Borel subgroup, and ${\bf T}$ be an $\text{Fr}$ -stable maximal torus contained in ${\bf B}$, and ${\bf U}=R_u({\bf B})$ be the ($\text{Fr}$ -stable) unipotent radical of ${\bf B}$. We identify ${\bf G}$ with ${\bf G}(\bar{\mathbb{F}}_q)$ and do likewise for the various subgroups of ${\bf G}$ such as ${\bf B}, {\bf T}, {\bf U}$ $\cdots$. We denote by $\Phi=\Phi({\bf G};{\bf T})$ the corresponding root system, and by $\Phi^+$ (resp. $\Phi^-$) the set of positive (resp. negative) roots determined by ${\bf B}$. Let $W=N_{\bf G}({\bf T})/{\bf T}$ be the corresponding Weyl group. We denote by $\Delta=\{\alpha_i\mid 1\le i\le r\}$ the set of simple roots. For each $\alpha\in\Phi$, let ${\bf U}_\alpha$ be the root subgroup corresponding to $\alpha$ and we fix an isomorphism $\varepsilon_\alpha: \bar{\mathbb{F}}_q\rightarrow{\bf U}_\alpha$ such that $t\varepsilon_\alpha(c)t^{-1}=\varepsilon_\alpha(\alpha(t)c)$ for any $t\in{\bf T}$ and $c\in\bar{\mathbb{F}}_q$. For $1\le i\le r$, we denote $\varepsilon_i:=\varepsilon_{\alpha_i}$ for short. For any $w\in W$, let ${\bf U}_w$ be the subgroup of ${\bf U}$ generated by all ${\bf U}_\alpha$  with $w(\alpha)\in\Phi^-$.  One can refer to \cite{Car} for  the structure theory of algebraic groups.

 For convenience of later discussion, we recall some results in \cite{Xi}. For each positive integer $k$, we denote $G_k$ for the $\f_{q^{k!}}$-points of $\bg$, and do likewise for $B_k,T_k,U_k,U_{w,k}$ etc. Now let  ${\bf H}$ be a group which has a sequence of subgroups $H_1, H_2, \dots, H_k, \dots $ such that  ${\bf H}=\displaystyle  \bigcup_{i=1}^{\infty} H_i$ and $H_i \subset H_j$ whenever $i<j$. Clearly, the groups ${\bf G}, {\bf B}, {\bf T}, {\bf U}$ satisfy this property. Let $\{M_i, \varphi_{ij}\}$ be a direct system of vector spaces over $\Bbbk$.  Assume that  $M_i$ is a $\Bbbk H_i$-module for each integer $i$ and the morphisms $\varphi_{ij}:$ $M_i\rightarrow M_j$ satisfy $\varphi_{ij}(hm)=h\varphi_{ij}(m)$ for any $h\in H_i$ and $m\in M_i$ whenever $i<j$.  Then the direct limit $\varinjlim M_i$ is a $\Bbbk {\bf H}$-module by \cite[Lemma 1.5]{Xi}. This is a useful method to construct the abstract representations of algebraic groups. Moreover, N.H. Xi showed that if each $M_i$ is an irreducible $\Bbbk H_i$-module, then the direct limit $\varinjlim M_i$ is an irreducible  $\Bbbk {\bf H}$-module.

 Let $\widehat{\bf T}$ be the set of characters of ${\bf T}$ over ${\Bbbk}$. Clearly, $W$ acts on $\widehat{\bf T}$ by $\theta\mapsto\theta^w$, where $\theta^w(t)=\theta(w^{-1}tw)$, $t\in{\bf T}$.  Each $\theta\in\widehat{\bf T}$ is regarded as a character of ${\bf B}$ by the homomorphism ${\bf B}\rightarrow{\bf T}$. Let ${\Bbbk}_\theta$ be the corresponding ${\bf B}$-module. We are interested in the induced module $\mathbb{M}(\theta)=\Bbbk{\bf G}\otimes_{\Bbbk{\bf B}}{\Bbbk}_\theta$. Fix a nonzero element ${\bf 1}_{\theta}$  in ${\Bbbk}_\theta$. We write $x{\bf 1}_{\theta}:=x\otimes{\bf 1}_{\theta}\in \mathbb{M}(\theta)$ for short. Using the Bruhat decomposition, we see that $$\mathbb{M}(\theta)= \sum_{w\in W}\Bbbk {\bf U}_{w^{-1}} \dot{w} {\bf 1}_{\theta},$$
 and thus it is not difficult to see that $\op{End}_{\bk\bg}(\m(\theta))=\bk$ since $\bk\one_\theta$ is the unique $\bb$-stable line in $\m(\theta)$ (see \cite[Proposition 2.2]{CD1}). Thus $\mathbb{M}(\theta)$ is an indecomposable $\Bbbk {\bf G}$-module.
 For each integer $i$, let $M_i(\theta)=\bk G_i\otimes_{\Bbbk B_i}{\Bbbk}_\theta$, which is a $\bk G_i$-module. Then $M_i(\theta)$ is identified with $\bk G_i\one_\theta$ and we have $\m(\theta)=\bigcup_iM_i(\theta)$.

For each $i \in I$, let ${\bf G}_i$ be the subgroup of $\bf G$ generated by ${\bf U}_{\alpha_i}, {\bf U}_{-\alpha_i}$ and we set ${\bf T}_i= {\bf T}\cap {\bf G}_i$. For $\theta\in\widehat{\bf T}$, define the subset $I(\theta)$ of $I$ by $$I(\theta)=\{i\in I \mid \theta| _{{\bf T}_i} \ \text {is trivial}\}.$$
For $J\subset I(\theta)$, let $W_J$ be the subgroup of $W$ generated by $s_i~(i\in J)$, and let ${\bf G}_J$ be the subgroup of $\bf G$ generated by ${\bf G}_i~(i\in J)$. We choose a representative $\dot{w}\in {\bf G}_J$ for each $w\in W_J$. Thus, the element $w{\bf 1}_\theta:=\dot{w}{\bf 1}_\theta$  $(w\in W_J)$ is well defined.

For $J\subset I(\theta)$, we set
$$\eta(\theta)_J=\sum_{w\in W_J}(-1)^{\ell(w)}w{\bf 1}_{\theta},$$
where $\ell(w)$ is the length of  $w\in W$.  Let $\mathbb{M}(\theta)_J=\Bbbk{\bf G}\eta(\theta)_J$ be the $\Bbbk {\bf G}$-module which is generated by $\eta(\theta)_J$.  The structure of  $\mathbb{M}(\theta)_J$ is very useful to our later discussion.  We state the results here. One can refer to \cite[Section 2]{CD1} for more details. Denote by $$W^J =\{x\in W\mid x~\op{has~minimal~length~in}~xW_J\}$$ and we have
$$\mathbb{M}(\theta)_J=\sum_{w\in W^J}\Bbbk {\bf U}_{w_Jw^{-1}}\dot{w}\eta(\theta)_J.$$
We give the formula of $\dot{s_i}  \varepsilon_i(x) \dot{w} \eta(\theta)_J$,  where $w\in W^J$ and $e\ne  \varepsilon_i(x)  \in {\bf U}_{\alpha_i} \subseteq  {\bf U}_{w_Jw^{-1}}$. The details are referred to \cite[Proposition 2.5]{CD1}.
For each $\varepsilon_i(x) \in {\bf U}_{\alpha_i}\backslash\{e\}$, we have
$$\dot{s_i}\varepsilon_i(x)\dot{s_i}^{-1}=f_i(x)\dot{s_i}h_i(x)g_i(x),$$
where $f_i(x),g_i(x) \in {\bf U}_{\alpha_i}\backslash\{e\}$, and $h_i(x)\in {\bf T}_i$ are uniquely determined.
Noting that  ${\bf U}_{\alpha_i} \subseteq  {\bf U}_{w_Jw^{-1}}$,   we get $s_i w w_J \leq ww_J$. By easy calculation, we have the following

\begin{proposition} \label{suwformual}
With the notations above,

\noindent (i) If $s_iw \leq w$, then
$$\dot{s_i} \varepsilon_i(x) \dot{w} \eta(\theta)_J=\theta^w(\dot{s_i}h_i(x)\dot{s_i})f_i(x)\dot{w}\eta(\theta)_J.$$

\noindent  (ii) If $w \leq  s_iw$ but $s_iww_J\leq ww_J$, then
$$\dot{s_i} \varepsilon_i(x)  \dot{w} \eta(\theta)_J =(f_i(x)-1)\dot{w}\eta(\theta)_J.$$
\end{proposition}

We define
$$E(\theta)_J=\mathbb{M}(\theta)_J/N(\theta)_J,$$
where $N(\theta)_J$ is the sum of all $\mathbb{M}(\theta)_K$ with $J\subsetneq K\subset I(\theta)$.
 If $\op{char}\bk\ne p$, then the composition factors of  $\mathbb{M}(\theta)$  are $E(\theta)_J$ $ (J\subset I(\theta))$,  with each of multiplicity one (see \cite[Theorem 3.1]{CD1}). According to \cite[Proposition 2.8]{CD1}, one has that  $E(\theta_1)_{K_1}$ is isomorphic to $E(\theta_2)_{K_2}$ as $\Bbbk {\bf G}$-modules if and only if $\theta_1=\theta_2$ and $K_1=K_2$.

In the paper  \cite{D1}, we introduce the principal representation category $\mathscr{O}({\bf G})$. The category $\mathscr{O}({\bf G})$  is  defined to be the full subcategory of $\Bbbk{\bf G}$-Mod such that any object $M$ in $\mathscr{O}({\bf G})$ is of finite length and its composition factors are $E(\theta)_J$ for some $\theta \in \widehat{\bf T}$ and $J\subset I(\theta)$.

\section{Extensions by finite-dimensional modules in $\mathscr{O}({\bf G})$}

 According to  a theorem of Borel and Tits (see \cite[Theorem 10.3 and Corollary 10.4]{BT}),  if ${\bf G}$ is a semisimple algebraic group and $\op{char}\bk\ne p$, then except the trivial representation, all other irreducible representations of $\Bbbk {\bf G}$ are infinite dimensional.  When  ${\bf G}$ is reductive,
  each finite-dimensional irreducible representation is one dimensional by  \cite[Corollary 6.6]{CD1}. So the finite-dimensional representations are very few among all the representations of ${\bf G}$.  The main result of this section is

\begin{theorem}\label{extbyfin}
Suppose that $p\ne \op{char}\bk\ge5$ or $\op{char}\bk=0$. If  $M,N\in\mathscr{O}(\bg)$ and $\dim N<\infty$, then $ \op{Ext}_{\bk\bg}^1(M, N)=0$.
\end{theorem}
\medskip

It is enough to prove Theorem \ref{extbyfin} in the case  $\dim N=1$ and $M=E(\theta)_J$. Replacing $N$ by $N\otimes\chi^{-1}$ and $E(\theta)_J$ by $E(\theta\otimes\chi^{-1}|_{\bt})_J$ when necessary ($\chi:\bg\rightarrow\bk^\times$ is the character corresponding to $N$), we can assume that $N={\rm tr}$, the trivial $\bg$-module. Thus, it is enough to show that
\begin{theorem}\label{extbyfin2}
Suppose that $p\ne \op{char}\bk\ge5$ or $\op{char}\bk=0$.  We have $\op{Ext}_{\bk\bg}^1(E(\theta)_J, {\rm tr})=0$ for any $\theta\in\widehat{\bt}$ and $J\subset I(\theta)$.
\end{theorem}

By the definition of $ E(\theta)_J$, we have the short exact sequence
$$ 0 \rightarrow N(\theta)_J \rightarrow  \mathbb{M}(\theta)_J  \rightarrow  E(\theta)_J \rightarrow 0,$$
which induces the long exact sequence
$$
\aligned
0 &\ \rightarrow  \text{Hom}_{\bk\bg}(E(\theta)_J, \text{tr})  \rightarrow  \text{Hom}_{\bk\bg}( \mathbb{M}(\theta)_J, \text{tr})  \rightarrow \text{Hom}_{\bk\bg}(N(\theta)_J, \text{tr}) \rightarrow  \\
&\  \rightarrow \op{Ext}_{\bk\bg}^1(E(\theta)_J, \text{tr})
\rightarrow  \op{Ext}_{\bk\bg}^1( \mathbb{M}(\theta)_J, \text{tr})
 \rightarrow\cdots.
\endaligned
$$
Since tr is not a composition factor of $N(\theta)_J$ by \cite[Theorem 3.1]{CD1}, we have $ \text{Hom}_{\bk\bg}(N(\theta)_J, \text{tr})=0$. So it is enough to show that $\op{Ext}_{\bk\bg}^1( \mathbb{M}(\theta)_J, \text{tr}) =0$, which is the main goal of the remaining part of this section.

\begin{lemma}\label{specialelement}
  Given a short exact sequence
\begin{equation}\label{exttrM}
0 \rightarrow \op{tr} \rightarrow  M \rightarrow  \mathbb{M}(\theta)_J \rightarrow 0,
\end{equation}
then there exists an element $\xi(\theta)_J\in M$ whose image in $ \mathbb{M}(\theta)_J $ is $\eta(\theta)_J$ with the following properties
$$t\xi(\theta)_J=\theta(t)\xi(\theta)_J, \ \forall t\in {\bf T};  \quad u\xi(\theta)_J= \xi(\theta)_J,\  \forall u\in {\bf U}'_{w_J}$$
and
$$ \dot{w} \xi(\theta)_J= (-1)^{\ell(w)} \xi(\theta)_J,  \quad  \forall w\in W_J.$$
\end{lemma}

\begin{proof}
Let $\xi(\theta)_J$ be an elements in $M$ whose image in $ \mathbb{M}(\theta)_J $ is $\eta(\theta)_J$. Then we have
$$ t \xi(\theta)_J=\theta(t) \xi(\theta)_J+ \varphi(t) {\bf 1}_{\text{tr}}, \quad \forall t\in {\bf T}.$$
If $\theta$ is trivial, then we get $\varphi(t_1t_2)= \varphi(t_1)+\varphi(t_2), \forall t_1, t_2\in {\bf T}.$ In the case that $\text{char}\ \Bbbk=l=0$, each $t\in {\bf T}$ is of finite order which is  denoted by $o(t)$ and thus we have
$$\varphi(e)=\varphi(t^{o(t)})=o(t)\varphi(t)= 0$$
and so $\varphi(t)=0, \forall t\in {\bf T}$. When $\text{char}\ \Bbbk=l>0$, there exists $z\in {\bf T}$ such that $z^{l}=t$  for any $t\in {\bf T}$. So we get
$\varphi(t)=l\varphi(z) =0$.

Now we consider that $\theta$ is non-trivial. It is easy to see that
$$ \varphi(xy)= \theta(y) \varphi(x)+ \varphi(y)=\theta(x) \varphi(y)+ \varphi(x), \  \forall x,y\in {\bf T}.$$
which implies that $$(\theta(y)-1)\varphi(x)= (\theta(x)-1)\varphi(y).$$
So $\displaystyle \frac{\varphi(x)}{\theta(x)-1}$ is a constant (denoted by $a$) for any $x$ with $\theta(x)\ne 1$, and hence we write $\varphi(x)=a(\theta(x)-1) $ for any $x\in {\bf T}$. Then we have
$$t(\xi(\theta)_J+a{\bf 1}_{\text{tr}})= \theta(t)(\xi(\theta)_J+a {\bf 1}_{\text{tr}}), \quad \forall t\in {\bf T}.  $$
So now we have gotten an element  which is still denoted by $\xi(\theta)_J$ such that $t\xi(\theta)_J=\theta(t)\xi(\theta)_J, \forall t\in {\bf T}$ and its image in $ \mathbb{M}(\theta)_J $ is $\eta(\theta)_J$.

Now we show that $\varepsilon_{\alpha} (x) \xi(\theta)_J= \xi(\theta)_J$ for any $\alpha \in \Phi- \Phi_J$. Suppose that
$$ \varepsilon_{\alpha} (x) \xi(\theta)_J= \xi(\theta)_J+ \psi_\alpha(x){\bf 1}_{\text{tr}}, \quad  \forall x\in \bar{\mathbb{F}}_q.$$
Then we have
$$\psi_\alpha(x+y)=\psi_\alpha(x)+ \psi_\alpha(y), \quad \forall  x, y\in \bar{\mathbb{F}}_q.$$
In particular, we have $$\psi_\alpha(e)= \psi_\alpha(px)=p\psi_\alpha(x)=0$$ and so $\psi_\alpha(x)=0$ for all $ x\in \bar{\mathbb{F}}_q$ since $p\ne 0$ in $\Bbbk$. In particular, we get $$u\xi(\theta)_J= \xi(\theta)_J,\quad  \forall u\in {\bf U}'_{w_J}.$$

For any $w\in W_J$,  we have
$$ \dot{w} \xi(\theta)_J= (-1)^{\ell(w)} \xi(\theta)_J+ \phi(w){\bf 1}_{\text{tr}}.$$
Thus for any simple reflection $s\in J$, we get  $\phi(ws)=\phi(s)-\phi(w)$. In particular
$$\phi((rs)^m)=m(\phi(s)-\phi(r)), \quad \forall s, r\in J.$$ Noting that $\text{char} \ \Bbbk =0$ or $\text{char}\ \Bbbk=l\geq 5$,   when $m$ is the order of $rs$ ($m=2,3,4,6$),  we have $\phi(s)$ is a constant (denoted by $b$) for any $s\in J$. Therefore we get
 $$ \dot{s} \xi(\theta)_J= -\xi(\theta)_J+ b{\bf 1}_{\text{tr}}.$$
Let $t\in {\bf T}$ act on both sides and then we have $b = \theta(t)b, \forall t\in {\bf T}$. If $\theta$ is non-trivial, then we get $b=0$ and thus $\phi(w)=0$ for all $w\in W_J$. If $\theta$ is trivial, denote by $c=\displaystyle -\frac{b}{2}$ and then
$$ \dot{w}(\xi(\theta)_J+ c{\bf 1}_{\text{tr}} )= (-1)^{\ell(w)}(\xi(\theta)_J+ c{\bf 1}_{\text{tr}} ),  \quad  \forall w\in W_J.$$
In this case,  $\xi(\theta)_J+ c{\bf 1}_{\text{tr}} $ is also a  ${\bf T}$-stable element. We also have
$$u(\xi(\theta)_J+ c{\bf 1}_{\text{tr}} )= \xi(\theta)_J+ c{\bf 1}_{\text{tr}},  \quad  \forall u\in {\bf U}'_{w_J}.$$
For convenience, this element is still denoted by $\xi(\theta)_J$.  Therefore the lemma is proved.
\end{proof}

Now we show that $\op{Ext}_{\bk\bg}^1( \mathbb{M}(\theta)_J, \text{tr}) =0$, which implies Theorem \ref{extbyfin2}. Given a short exact sequence  as (\ref{exttrM}), let $\xi(\theta)_J$ be the element which satisfies the conditions in Lemma \ref{specialelement}.  It remains to show that  the subspace $M'$ of $M$ spanned by the following set
$$\{ u \dot{w} \xi(\theta)_J \mid w\in W^J, u\in{\bf U}_{w_Jw^{-1}} \}$$
is a submodule of $M$ which is isomorphic to $\mathbb{M}(\theta)_J $. Thus the short exact sequence  (\ref{exttrM})  is splitting and we have $\text{Ext}_{\bk\bg}^1( \mathbb{M}(\theta)_J, \text{tr}) =0$. To prove that $M'$ is a submodule of $M$, it is enough to show that
$\dot{s_i}  \varepsilon_i(x) \dot{w} \xi(\theta)_J\in M',$ where $w\in W^J$ and $e\ne  \varepsilon_i(x)  \in {\bf U}_{\alpha_i} \subseteq  {\bf U}_{w_Jw^{-1}}$.

If $s_iw \leq w$, using Proposition \ref{suwformual}, then we have
$$\dot{s_i} \varepsilon_i(1) \dot{w} \xi(\theta)_J=\theta^w(\dot{s_i}h_i(1)\dot{s_i})f_i(1)\dot{w}\xi(\theta)_J+ \kappa_i {\bf 1}_{\text{tr}}.$$
Using ${\bf T}$ acting on both sides and then for any $ x\in \bar{\mathbb{F}}^*_q$,   we have
$$\dot{s_i} \varepsilon_i(x) \dot{w} \xi(\theta)_J=\theta^w(\dot{s_i}h_i(x)\dot{s_i})f_i(x)\dot{w}\xi(\theta)_J+ \kappa_i {\bf 1}_{\text{tr}},$$
and symmetrically, we get
$$\dot{s_i} f_i(x) \dot{w} \xi(\theta)_J=\theta^w(\dot{s_i}h_i(x)\dot{s_i})^{-1}\dot{s}^2\varepsilon_i(x)\dot{w}\xi(\theta)_J+ \kappa_i {\bf 1}_{\text{tr}}.$$
Combining these two equations we get $\theta^w(\dot{s_i}h_i(x)\dot{s_i}) \kappa_i +\kappa_i=0$ for any $x\in \bar{\mathbb{F}}^*_q$.
Since $\theta^w(y)\ne -1$ for some $y\in{\bf T}_i$, we have $\kappa_i=0$.

If $w \leq  s_iw$ but $s_iww_J\leq ww_J$,  using Proposition \ref{suwformual} and by the same discussion  we have
$$\dot{s_i} \varepsilon_i(x)  \dot{w} \xi(\theta)_J =(f_i(x)-1)\dot{w}\xi(\theta)_J+ \vartheta_i {\bf 1}_{\text{tr}}, \  \forall x\in \bar{\mathbb{F}}^*_q.$$
Let $\dot{s_i}$ act on both sides and we get $2  \vartheta_i=0$ and so $ \vartheta_i=0$. Therefore we have proved that $M'$ is a submodule of $M$ and the sequence (\ref{exttrM}) splits and Theorem \ref{extbyfin2} is proved.

\section{Extension of naive induced modules}

The study  of  induced modules  $\m(\theta)$ is motivated by the Verma modules. The extensions of Verma modules are very important in the study of BGG category $\mathcal{O}$.  So we are also interested in the extension of induced modules.  The main result of this section is

\begin{theorem} \label{mainthm}
 Let $\lambda,\mu\in\widehat{\bf T}$. Then $\op{Ext}_{\bk\bg}^1(\m(\lambda),\m(\mu))=0$ if and only if $\lambda|_{\bf C}\ne\mu|_{\bf C}$.
\end{theorem}
The above theorem is true for any field $\Bbbk$. In \cite{D2}, J.B. Dong  introduced a full-subcategory $\mathscr{X}(\bf G)$ of $\Bbbk {\bf G}$-Mod when $\Bbbk =\mathbb{C}$. In this category $\mathscr{X}(\bf G)$,   we see that  $\op{Ext}_{\mathscr{X}(\bf G)}^1(\m(\lambda),\m(\mu))=0$ whenever $\lambda \ne \mu$ by \cite[Proposition 6.8]{D2}. However Theorem \ref{mainthm} shows that the extensions in $\Bbbk {\bf G}$-Mod are very complicated.

 Since the functor $\Bbbk{\bf G}\otimes_{\Bbbk{\bf B}}-$ is exact and takes projective $\Bbbk{\bf B}$-modules to projective $\Bbbk{\bf G}$-modules, we have $$\op{Ext}_{\bk\bg}^1(\m(\lambda),\m(\mu))\simeq\op{Ext}_{\bk\bb}^1(\Bbbk_\lambda,\m(\mu)).$$ Thus Theorem \ref{mainthm} is equivalent to the following

 \begin{theorem}\label{thm1}
One has   $\op{Ext}_{\bk\bb}^1(\Bbbk_\lambda,\m(\mu))=0$ if and only if $\lambda|_{\bf C}\ne\mu|_{\bf C}$.
 \end{theorem}

 This section is devoted  to prove Theorem \ref{thm1}. Firstly,  we prove the ``if" part of Theorem \ref{thm1}.

\begin{proposition} \label{ifpart}
If $\lambda|_{\bf C}\ne\mu|_{\bf C}$, then $\op{Ext}_{\bk\bb}^1(\Bbbk_\lambda,\m(\mu))=0$.
\end{proposition}

\begin{proof}
Assume that $\lambda|_{\bf C}\ne\mu|_{\bf C}$, then there is a $c_0\in\mathbb{\bf C}$ such that $\lambda(c_0)\ne\mu(c_0)$. Let
\begin{equation}\label{extb}
0\rightarrow \m(\mu)\rightarrow E\rightarrow\Bbbk_\lambda\rightarrow0
\end{equation}
be an exact sequence of $\Bbbk{\bf B}$-modules.
Let $\one_\lambda'$ be the pre-image of $\one_\lambda$ in $E$. Then one has $c_0\one_\la'=\la(c_0)\one_\la'+m_0$ for some $m_0\in\m(\mu)$. For any $c\in\bc$, the map $m\mapsto cm$ ($m\in\m(\mu)$) is in $\op{End}_{\bk\bg}(\m(\mu))$. Note that  $\op{End}_{\bk\bg}(\m(\mu))=\bk$ and thus the action of $\bc$ on $\m(\mu)$ is a scalar, and hence $cm=\mu(c)m$ for any $m\in\m(\mu)$ since $c\one_\mu=\mu(c)\one_\mu$.

Let $\alpha=\one_\la'+(\la(c_0)-\mu(c_0))^{-1}m_0\in E$. Then
$$
\aligned
c_0\alpha &\ =c_0\one_\la'+(\la(c_0)-\mu(c_0))^{-1}c_0m_0\\
&\ =\la(c_0)\one_\la'+m_0+\mu(c_0)(\la(c_0)-\mu(c_0))^{-1}m_0\\
&\ =\la(c_0)(\one_\la'+(\la(c_0)-\mu(c_0))^{-1}m_0)\\
&\ =\la(c_0)\alpha.
\endaligned
$$
Let $E_{c_0,\la(c_0)}=\{v\in E\mid c_0v=\la(c_0)v\}$, which is non-zero since $\alpha \in E_{c_0,\la(c_0)}$.  Noting that $\la(c_0)\ne\mu(c_0)$, we have $E_{c_0,\la(c_0)}\cap\m(\mu)=0$, and hence $\dim E_{c_0,\la(c_0)}=1$. It is clear that $E_{c_0,\la(c_0)}$ is a $\bk\bb$-module since $c_0$ commutes with each element in $\bb$. Thus, we have the split exact sequence
$$0\rightarrow \m(\mu)\rightarrow E\rightarrow E_{c_0,\la(c_0)}\rightarrow0$$
of $\bk\bb$-modules, and hence $E_{c_0,\la(c_0)}\simeq E/\m(\mu)\simeq\bk_\la$. This implies that  (\ref{extb}) is split and thus  the ``if" part of
Theorem  \ref{thm1} is proved.
\end{proof}

In the following, we prove the ``only if" part of Theorem \ref{thm1}. From now on we assume that $\la|_{\bc}=\mu|_{\bc}$. It is enough to construct a non-split extension $F$ of $\bk_\la$ and $\m(\mu)$. We make some preliminaries before constructing $F$.

There is an natural group homomorphism
$$\pi:~\bu\rightarrow\bu_{\alpha_1}\times\bu_{\alpha_2}\times\cdots\times\bu_{\alpha_r}$$
whose kernel is the subgroup of $\bu$ generated by $\bu_\alpha$ with $\alpha\in\Phi^+\backslash\Delta$. It is clear that $\pi$ maps $U_i$ to $U_{\alpha_1,i}\times\cdots\times U_{\alpha_r,i}$ and compatible with the conjugation of $\bt$.

For each finite subset $H$ of $\bg$, set $\underline{H}=\sum_{h\in H}h\in\bk\bg$. We fix an element $a_i\in\f_{q^{(i+1)!}}-\f_{q^{i!}}$ and define $$u_i=\varepsilon_1(a_i)\varepsilon_2(a_i)\cdots\varepsilon_r(a_i)\in U_{i+1}$$ for each $i$. Let $N_i$ be the $\bk T_i$-module generated by $\underline{U_i}u_iw_0\one_\mu$. Then it is clear that $N_i$ is a $\bk T_i$-submodule of $M_{i+1}(\mu)$. Abbreviating $t\in T_i/T_i\cap\bc$ for ``$t$ runs over a complete set of representatives of cosets in $T_i/T_i\cap\bc$", then we have the following
\begin{lemma}\label{basis}
The set $\{\underline{U_i}tu_it^{-1}w_0\one_\mu\mid t\in T_i/T_i\cap\bc\}$ forms a basis of $N_i$.
\end{lemma}

\begin{proof}
It suffices to prove that $U_itu_it^{-1}$ $(t\in T_i/T_i\cap\bc)$ are distinct cosets in $U_i\backslash U_{i+1}$.
In fact, if $t_1,t_2\in T_i$ such that  $U_it_1u_it_1^{-1}=U_it_2u_it_2^{-1}$, then $t_1u_it_1^{-1}t_2u_i^{-1}t_2^{-1}\in U_i$, and hence we have $$\pi(t_1u_it_1^{-1}t_2u_i^{-1}t_2^{-1})\in U_{\alpha_1,i}\times\cdots\times U_{\alpha_r,i}.$$ On the other hand, we have
$$
\aligned
\pi(t_1u_it_1^{-1}t_2u_i^{-1}t_2^{-1}) &\ =\pi(t_1u_it_1^{-1})\pi(t_2u_i^{-1}t_2^{-1})\\
&\ =(\varepsilon_1(a_i\alpha_1(t_1)),\cdots,\varepsilon_r(a_i\alpha_r(t_1)))\\
&\ \quad\cdot(\varepsilon_1(-a_i\alpha_1(t_2)),\cdots,\varepsilon_r(-a_i\alpha_r(t_2)))\\
&\ =(\varepsilon_1(a_i(\alpha_1(t_1)-\alpha_1(t_2))),\cdots,\varepsilon_r(a_i(\alpha_r(t_1)-\alpha_r(t_2)))).
\endaligned
$$
It follows that $a_i(\alpha_k(t_1)-\alpha_k(t_2))\in \f_{p^{i!}}$ $(1\le k\le r)$. Since $a_i\not\in\f_{q^{i!}}$ and $\alpha_k(t_1)-\alpha_k(t_2)\in\f_{q^{i!}}$, we have $\alpha_k(t_1)=\alpha_k(t_2)$ $(1\le k\le r)$. Thus, we obtain $t_1t_2^{-1}\in\bigcap_{1\le k\le r}\op{Ker}\alpha_k=\bc$, and hence $t_1t_2^{-1}\in T_i\cap\bc$.

Conversely, if $t_1t_2^{-1}\in T_i\cap\bc$, then it is clear that $U_it_1u_it_1^{-1}=U_it_2u_it_2^{-1}$. This completes the proof.
\end{proof}

Let $\nu=(\mu^{w_0})^{-1}\la\in\widehat{\bt}$. Since $\bc$ is $W$-invariant and $\la|_\bc=\mu|_\bc$, we have $\nu|_\bc$ is trivial. It follows that the element
$$\eta_i=\sum_{t\in T_i/T_i\cap\bc}\nu(t)^{-1}\underline{U_i}tu_it^{-1}w_0\one_\mu\in M_{i+1}(\mu)$$
is well defined and $\eta_i\ne0$  by Lemma \ref{basis}. We have $u\eta_i=\eta_i$ for any $u\in U_i$, and
$$
\aligned
t_0\eta_i &\ =\mu^{w_0}(t_0)\sum_{t\in T_i/T_i\cap\bc}\nu(t)^{-1}\underline{U_i}t_0tu_i(t_0t)^{-1}w_0\one_\mu\\
&\ =\mu^{w_0}(t_0)\sum_{t\in T_i/T_i\cap\bc}\nu(tt_0^{-1})^{-1}\underline{U_i}tu_it^{-1}w_0\one_\mu\\
&\ =\mu^{w_0}(t_0)\nu(t_0)\eta_i\\
&\ =\la(t_0)\eta_i
\endaligned
$$
for any $t_0\in T_i$. In other words, the map $\one_\la\mapsto \eta_i$ is in $\op{Hom}_{B_i}(\bk_\la,M_{i+1}(\mu))$.

\smallskip

We start to construct a $\Bbbk {\bf B}$-module $F$. Firstly, define the $\bk B_i$-module $F_i=\bk_\la\oplus M_i(\mu)$. The map $f_i:~F_i\rightarrow F_{i+1}$  defined by
 $$ f_i(a\one_\la+m)=a\one_\la+a\eta_i+m, \quad a\in\bk, m\in M_i(\mu)$$ is in $\op{Hom}_{B_i}(F_i,F_{i+1})$ by the above discussion. For any $i<j$, we denote by $f_{i,j}:~F_i\rightarrow F_j$ the composition $f_{j-1}\circ\cdots f_{i+1}\circ f_i$ which is clearly in $\op{Hom}_{B_i}(F_i,F_j)$. Then $\{F_i,f_{i,j}\}$ is a direct system of vector spaces. In addition, $f_i(a\one_\la+m)=0$ implies that $a\one_\la=a\eta_i+m=0$, and hence $a=0$ and $m=0$. Thus, $f_i$ is injective and hence all $f_{i,j}$ are injective. Let $F$ be the direct limit of the direct system $\{F_i,f_{i,j}\}$.

Suppose that $v,v'\in F$ are represented by $v_i\in F_i$ and $v'_j\in F_j$ respectively. Then $v=v'$ if and only if $f_{i,k}(v_i)=f_{j,k}(v'_j)$ for some $k$. By the injectivity mentioned above, this is equivalent to
 $$
 \begin{cases}
 f_{i,j}(v_i)=v_j' &\ {\rm if}~i<j;\\
 v_i=v_i'          &\ {\rm if}~i=j;\\
 f_{j,i}(v_j')=v_i &\ {\rm if}~i>j.
 \end{cases}
 $$
 According to \cite[Lemma 1.5]{Xi}, $F$ is a $\bk\bb$-module via the following way. For any $b\in \bb$ and $v\in F$, suppose that  $v$ is represented by $v_i\in F_i$. We choose $k$ large enough such that $k\ge i$ and $b\in B_k$. Define $bv\in F$ to be the element represented by $bf_{ik}(v_i)\in F_k$.

\begin{lemma}\label{22}
There is an exact sequence
\begin{equation}\label{extb2}
0\rightarrow \m(\mu)\rightarrow F\rightarrow\Bbbk_\lambda\rightarrow0
\end{equation}
of $\bk\bb$-modules.
\end{lemma}
\begin{proof}
 We let
 $$M'=\{v\in F\mid v~{\rm is~represented~by~an~element~in}~M_k(\mu)~{\rm for~some}~k\}.$$ Since the restriction of $f_i$ on $M_i(\mu)$ is simply the usual inclusion $M_i(\mu)\rightarrow M_{i+1}(\mu)$, we have $M'\simeq\m(\mu)$ as $\bk\bb$-modules. Let $v,v'$ be the  nonzero element in $F$ and there is an integer $k$ such that $v,v'$ are represented by $v_k,v_k'\in F_k$ respectively. There is a $c\in\bk$ such that $v_k'-cv_k\in M_k(\mu)$ which implies $v'-cv\in M'$, and hence $\dim F/M'=1$. For any $b\in\bb$ and $v\in F$, there is an integer  $l$ such that $b\in B_l$ and $v$ is represented by $v_l\in F_l$. We have $bv_l-\la(b)v_l\in M_l(\mu)$ which implies $bv-\la(b)v\in M'$. Thus, we have $F/M'\simeq\bk_\la$. This completes the proof.
 \end{proof}

It remains to show that the following
 \begin{lemma}\label{23}
 The $\Bbbk {\bf B}$-module $F$ makes  the  sequence (\ref{extb2})  non-split.
 \end{lemma}
 \begin{proof}
 Suppose that $F\simeq\bk_\la\oplus\m(\mu)$, then there is a $0\ne v\in F^{\bu}$ with $v\not\in\m(\mu)$. Assume that $v$ is represented by $v_k=a\one_\la+m\in F_k$ $(a\in\bk,m\in M_k(\mu))$, and hence is also represented by $f_k(v_k)=a\one_\la+a\eta_k+m\in F_{k+1}$. We have $a\ne0$ since $v\not\in\m(\mu)$, and $f_k(v_k)\in F_{k+1}^{U_{k+1}}$ since $v\in F^\bu$. It follows that $a\eta_k+m\in F_{k+1}^{U_{k+1}}$. Let $m=\sum_{u\in U_k}c_uuw_0\one_\mu+m'$ $(c_u\in\bk)$, where $m'$ is a linear combination of $u'w\one_\mu$ with $w\ne w_0$ and $u'\in U_{w^{-1},k}$. Then we have $a\eta_k+\sum_{u\in U_k}c_uuw_0\one_\mu\in F_{k+1}^{U_{k+1}}$, and hence
\begin{equation}\label{inuk+1}
 a\eta_k+\sum_{u\in U_k}c_uuw_0\one_\mu\in\bk\underline{U_{k+1}}w_0\one_\mu.
\end{equation}
We claim that
\begin{equation}\label{clm}
\bigcup_{t\in T_k/T_k\cap\bc}U_ktu_kt^{-1}\ne U_{k+1}.
\end{equation}
In fact, $\bg=\bg'{\bf Z}$ is a product of a semisimple group $\bg'$ and a central torus ${\bf Z}$, and we get a corresponding decomposition $\bt={\bf S}{\bf Z}$, where ${\bf S}=\bt\cap\bg'$ is a maximal torus of $\bg'$. It follows that $|T_k/T_k\cap\bc|\le|{\bf S}^{{\text{Fr}}^{k!}}|\le(q^{k!}-1)^r$.
Since $|U_{k} \backslash U_{k+1}|=q^{(k+1)!|\Phi^+|}/q^{k!|\Phi^+|}=q^{k\cdot k!|\Phi^+|}$, we have
$$\left|\bigcup_{t\in T_k/T_k\cap\bc}U_ktu_kt^{-1}\right|=|U_k||T_k/T_k\cap\bc|\le(q^{k!}-1)^r|U_k|< q^{k\cdot k!|\Phi^+|}|U_k|=|U_{k+1}|$$
which implies the claim. Combining $a\ne0$ and (\ref{clm}) we see that
$$a\eta_k+\sum_{u\in U_k}c_uuw_0\one_\mu\not\in\bk\underline{U_{k+1}}w_0\one_\mu$$
which contradicts to (\ref{inuk+1}).
\end{proof}
Combining Lemma \ref{22} and \ref{23} we obtain Theorem \ref{thm1} and thus Theorem \ref{mainthm} is proved.

\section{Extensions of simple modules in $\mathscr{O}({\bf G})$ for $\bg=SL_2(\bar{\mathbb{F}}_q)$}

In this section, we assume that ${\bf G}= SL_2(\bar{\mathbb{F}}_q)$ and $p \neq \op{char}\Bbbk\geq 3$ or $\op{char}\bk=0$.  We consider the extensions of simple modules  in $\mathscr{O}({\bf G})$. For convenience of later discussion, we denote $$ h(t)=\begin{pmatrix}t&0\\0&t^{-1}\end{pmatrix}, \quad \varepsilon(a)=\begin{pmatrix}1&a\\0&1\end{pmatrix}, \quad s=\begin{pmatrix}0&-1\\1&0\end{pmatrix},$$
where $t\in \bar{\mathbb{F}}_q^\times$ and $a\in \bar{\mathbb{F}}_q$.

When $\theta$ is nontrivial, the induced module $\m(\theta)$ is already irreducible.  The module $\m(\text{tr})$ has the Steinberg module $\text{St}$ which is the $\bk\bg$-module generated by $\eta=(1-s)\one_{\rm tr}$ as its submodule and the quotient is the trivial module.  Now we denote
$$\text{Irr}({\bf G, \bf T}) =\{ \Bbbk_{\text{tr}}, \text{St}, \mathbb{M}(\theta)\mid \theta \in \widehat{{\bf T}} \ \text{is nontrivial}\}.$$
Moreover, the paper \cite{D3} has classified all the simple modules with ${\bf T}$-stable lines when $\Bbbk$ is  an algebraically closed field of characteristic zero. They are exactly the set $\text{Irr}({\bf G, \bf T})$. In this section, we determine the necessary and sufficient condition of the vanishing of the extensions between two simple modules in $\text{Irr}({\bf G, \bf T})$.

 The maximal torus  $\bt$ is identified with $\bar{\mathbb{F}}_q^\times$ and the center $\bc$ is identified with $\{\pm1\}$.   An easy calculation shows that
\begin{equation}\label{sus}
s\varepsilon(a)s=
\varepsilon(-a^{-1})sh(a)\varepsilon(-a^{-1}),\quad a\in\bar{\mathbb{F}}_q^\times.
\end{equation}

Firstly, by Theorem \ref{extbyfin}, we have $\op{Ext}_{\bk\bg}^{1}(M, \text{tr})=0$ for any $M\in \text{Irr}({\bf G, \bf T})$.
When $\theta$ is  nontrivial, using the short exact sequence
\begin{equation} \label{Sttrexact}
0 \rightarrow \text{St} \rightarrow \mathbb{M}(\text{tr}) \rightarrow \text{tr} \rightarrow 0,
\end{equation}
we get a long exact sequence
$$
\aligned
\cdots \rightarrow &\ \text{Hom}_{\bk\bg}(\mathbb{M}(\theta), \text{tr}) \rightarrow \op{Ext}_{\bk\bg}^{1}(\mathbb{M}(\theta), \text{St})\rightarrow
\op{Ext}_{\bk\bg}^{1}(\mathbb{M}(\theta),  \mathbb{M}(\text{tr})) \rightarrow  \\   &\  \rightarrow \op{Ext}_{\bk\bg}^{1}(\mathbb{M}(\theta),  \text{tr})  \rightarrow \cdots.
\endaligned
$$
Noting  that $ \text{Hom}_{\bk\bg}(\mathbb{M}(\theta), \text{tr})=  \op{Ext}_{\bk\bg}^{1}(\mathbb{M}(\theta),  \text{tr}) =0$,  we have $$\op{Ext}_{\bk\bg}^{1}(\mathbb{M}(\theta), \text{St}) \cong \op{Ext}_{\bk\bg}^{1}(\mathbb{M}(\theta), \mathbb{M}(\text{tr})).$$
Therefore $ \op{Ext}_{\bk\bg}^{1}(\mathbb{M}(\theta), \text{St})=0$  if and only if  $\theta|_{\bc}\ne {\rm tr}$ by Theorem \ref{mainthm}. It is obvious that $  \op{Ext}_{\bk\bg}^{1}(\text{tr}, \text{St})\ne 0$ since $\mathbb{M}(\text{tr})$ is already  a  nontrivial extension of  the trivial module $\text{tr}$ by
the Steinberg module $\text{St}$.
So we just need to consider the following extensions
$$  \op{Ext}_{\bk\bg}^{1}(\text{St}, \text{St}), \quad \op{Ext}_{\bk\bg}^{1}(\text{tr}, \mathbb{M}(\theta)),  \quad  \op{Ext}_{\bk\bg}^{1}(\text{St}, \mathbb{M}(\theta)).$$
In the following, we show that
 $\op{Ext}_{\bk\bg}^{1}(\text{tr}, \mathbb{M}(\theta))=0$ if and only if  $\theta|_{\bc}\ne {\rm tr}$,   and $\op{Ext}_{\bk\bg}^{1}(\text{St}, \mathbb{M}(\theta))=0$ if and only if  $\theta|_{\bc}\ne {\rm tr}$. In particular,  we get $\op{Ext}_{\bk\bg}^{1}(\text{St},  \mathbb{M}(\text{tr}))\ne 0$.
 Thus using the short exact sequence (\ref{Sttrexact}), we have a long exact sequence
 $$\cdots \rightarrow \text{Hom}_{\bk\bg}(\text{St}, \text{tr}) \rightarrow \op{Ext}_{\bk\bg}^{1}(\text{St}, \text{St})\rightarrow \op{Ext}_{\bk\bg}^{1}(\text{St},  \mathbb{M}(\text{tr}))  \rightarrow \cdots$$
which implies that $\op{Ext}_{\bk\bg}^{1}(\text{St}, \text{St}) \ne 0$.

\begin{proposition}
  If $\theta|_{\bc}\ne {\rm tr}$, then we have  $\op{Ext}_{\bk\bg}^{1}(\op{tr}, \mathbb{M}(\theta))=0$  and $\op{Ext}_{\bk\bg}^{1}(\op{St}, \mathbb{M}(\theta))=0$.
\end{proposition}

\begin{proof} Using the short exact sequence (\ref{Sttrexact}), we get a long exact sequence
$$\cdots \rightarrow \text{Hom}_{\bk\bg}(\text{St}, \mathbb{M}(\theta)) \rightarrow \op{Ext}_{\bk\bg}^{1}(\text{tr}, \mathbb{M}(\theta))\rightarrow \op{Ext}_{\bk\bg}^{1}(\mathbb{M}(\text{tr}), \mathbb{M}(\theta))  \rightarrow \cdots.$$
If $\theta|_{\bc}\ne {\rm tr}$, then $\op{Ext}_{\bk\bg}^{1}(\mathbb{M}(\text{tr}), \mathbb{M}(\theta)) =0$. Noting that $\text{Hom}(\text{St}, \mathbb{M}(\theta))=0$, we get $\op{Ext}_{\bk\bg}^{1}(\text{tr}, \mathbb{M}(\theta))=0$.

Now we consider the following short exact sequence
$$0 \rightarrow\mathbb{M}(\theta)  \xrightarrow{\phi} M \xrightarrow{\psi} \text{St} \rightarrow 0.$$
Since $\theta|_{\bc}\ne {\rm tr}$, we get $\op{Ext}_{\bk\bt}^{1}(\text{tr}, \mathbb{M}(\theta))=0$ by the similar arguments  in the proof of Proposition \ref{ifpart}. Thus there exists $\xi \in M$ such that $t\xi =\xi$ for any $t\in {\bf T}$ and $\psi(\xi)=\eta$.  Assume that
$s\xi=-\xi +z$, where $z\in \mathbb{M}(\theta)$. Let $t\in {\bf T}$ act on both sides. Then we get $tz=z$ for any $t\in {\bf T}$, which implies that $z=0$.  Now we consider
$$s \varepsilon(x) \xi=( \varepsilon(-x^{-1})-1)\xi + v_x,$$
where $v_x\in \mathbb{M}(\theta)$. However let $h(-1)\in {\bf T}$ act on both side, we have $h(-1) v_x= v_x$, which implies that $v_x=0$ since $\theta|_{\bc}\ne {\rm tr}$.
Thus $\Bbbk {\bf U}\xi$ is a submodule of $M$, which is isomorphic to the Steinberg module $\text{St}$. So we get $\op{Ext}_{\bk\bg}^{1}(\text{St}, \mathbb{M}(\theta))=0$.
\end{proof}

Since $\theta|_{\bc}={\rm tr}$, $\theta$ induces a character of $\bt/\{\pm1\}$. Thus one can assume that $\bg=PGL_2(\bar{\mathbb{F}}_q)=SL_2(\bar{\mathbb{F}}_q)/\{\pm1\}$ in which one identifies $T_i$ with $\f_{q^{i!}}^\times/\{\pm1\}$, and we abbreviate $t\in T_i$ for ``choosing a representative $t$ in $\f_{q^{i!}}^\times$" and $\theta(t)$ for $\theta(h(t))$. We keep this assumption in the sequel.

\begin{proposition}\label{main2}
If $\theta|_{\bc}={\rm tr}$, then $\op{Ext}_{\bk\bg}^{1}({\rm tr},\m(\theta))\ne0$
\end{proposition}

 To prove this proposition, we will construct a nonsplit exact sequence
 \begin{equation}\label{exth}
0\rightarrow\m(\theta)\rightarrow H\rightarrow{\rm tr}\rightarrow0
 \end{equation}
of $\bk\bg$-modules.  Let $H_i={\text{tr}}\oplus M_i(\theta)$ which is a $\bk G_i$-module. For each $i$, we fix an element  $b_i\in\f_{q^{(i+1)!}}-\f_{q^{2i!}}$ and set $$\xi_i=\underline{G_i}\varepsilon(b_i)\sn\in M_{i+1}(\theta).$$
\begin{lemma}\label{ne02}
We have $\xi_i\ne0$ and $g\xi_i=\xi_i$ for any $g\in G_i$.
\end{lemma}
\begin{proof}
The second equality is trivial and we just prove the first one. By the Bruhat decomposition, we have
$$
\xi_i=\underline{B_i}\varepsilon(b_i)\sn+\underline{U_i}s\underline{B_i}\varepsilon(b_i)\sn.
$$
The first part is \begin{equation}\label{1}
\underline{B_i}\varepsilon(b_i)\sn=\underline{U_i}\cdot\underline{T_i}\varepsilon(b_i)\sn
=\underline{U_i}\sum_{t\in T_i}\theta(t)^{-1}\varepsilon(b_it^2)\sn,
\end{equation}
and the second part is
\begin{equation}\label{2}
\aligned
s\underline{B_i}\varepsilon(b_i)\sn &\ =s\underline{U_i}\sum_{t\in T_i}\theta(t)^{-1}\varepsilon(b_it^2)\sn\\
&\ =\sum_{{t\in T_i}\atop{a\in\f_{q^{i!}}}}\theta(t)^{-1}s\varepsilon(b_it^2+a)\sn\\
&\ =\sum_{{t\in T_i}\atop{a\in\f_{q^{i!}}}}\theta(t)^{-1}\theta(b_it^2+a)\varepsilon(-(b_it^2+a)^{-1})\sn.
\endaligned
\end{equation}
Combining (\ref{1})  and  (\ref{2}),  we have
$$\xi_i=\underline{U_i}\left[\sum_{t\in T_i}\theta(t)^{-1}\varepsilon(b_it^2)\sn+\sum_{{t\in T_i}\atop{a\in\f_{q^{i!}}}}\theta(t)^{-1}\theta(b_it^2+a)\varepsilon(-(b_it^2+a)^{-1})\sn\right].$$
We claim that $$\{U_i\varepsilon(b_it^2) \mid t\in T_i\} \cup  \{U_i\varepsilon(-(b_it^2+a)^{-1}) \mid t\in T_i, a\in\f_{q^{i!}}\}$$ are distinct cosets in $U_i\backslash U_{i+1}$. In fact, let
$$x,y\in\{b_it^2\mid t\in T_i\}\cup\{-(b_it^2+a)^{-1}\mid t\in T_i, a\in\f_{q^{i!}}\}$$
and $x\ne y$. Suppose that $\varepsilon(x)\varepsilon(y)^{-1}\in U_i$. Then $x-y\in\f_{q^{i!}}$, which is equivalent to $b_i$ is a root of certain quadratic (or linear) equation with coefficients in $\f_{q^{i!}}$ by simple calculation. This contradicts to the setting of  $b_i\in\f_{q^{(i+1)!}}-\f_{q^{2i!}}$, and the claim holds. In particular, we have $\xi_i\ne0$.
\end{proof}
Clearly, the map $h_i: H_i\rightarrow H_{i+1}$ defined by $$h_i(a {\bf 1}_{\text{tr}}+m)=a {\bf 1}_{\text{tr}} +a\xi_i+m, \quad  a\in\bk, m\in M_i(\theta)$$ is an injective homomorphism in $\op{Hom}_{G_i}(H_i,H_{i+1})$ by Lemma \ref{ne02}. For any $i<j$, we denote by $h_{i,j}:~H_i\rightarrow H_j$ the composition $h_{j-1}\circ\cdots\circ h_{i+1}\circ h_i$ which is clearly in $\op{Hom}_{G_i}(H_i,H_j)$. Then $\{H_i,h_{i,j}\}$ is a direct system of vector spaces. Let $H$ be the direct limit of $\{H_i,h_{i,j}\}$, which is a $\Bbbk {\bf G}$-module.  The same arguments as in the proof of Lemma \ref{22} shows that $H$ satisfies the exact sequence (\ref{exth}).

\begin{lemma}\label{noHG}
There is no $v\in H-\m(\theta)$ such that $v\in H^{\bg}$.
\end{lemma}
\begin{proof}
Suppose $v\in H-\m(\theta)$ and $v\in H^{\bg}$, and assume that $v$ is represented by $v_i\in H_i$. Then $h_i(v_i)\in H_{i+1}^{G_{i+1}}$ since $v\in H^{\bg}$. Let $v_i=a {\bf 1}_{\text{tr}} +m$, where $a\in\bk, m\in\m(\theta)$. Then $a\ne0$ since we have $v\in H-\m(\theta)$. Noting that  $$H_{i+1}^{G_{i+1}} \subset H_{i+1}^{U_{i+1}}=\bk {\bf 1}_{\text{tr}} +\bk\one_\theta+\bk\underline{U_{i+1}}\sn$$ and $a\ne0,\xi_i\ne0$ (by Lemma \ref{ne02}), then we get $h_i(v_i)\not\in\bk {\bf 1}_{\text{tr}}+\bk\one_\theta$, i.e.,
\begin{equation}\label{35}
 a {\bf 1}_{\text{tr}} +a\xi_i+m=a {\bf 1}_{\text{tr}}+x\one_\theta+y\underline{U_{i+1}}\sn,
\end{equation}
where $y\ne 0$.  On the other hand, we have
\begin{equation}\label{36}
\bigcup_{t\in T_i}U_i\varepsilon(b_it^2)\cup\bigcup_{{t\in T_i}\atop{a\in\f_{q^{i!}}}}U_i\varepsilon(b_it^2+a)\ne U_{i+1}
\end{equation}
since the cardinality of the left side of (\ref{36}) is $$|U_i|\left(\displaystyle\frac{q^{i!}-1}{2}+q^{i!}\displaystyle\frac{q^{i!}-1}{2}\right)=\displaystyle\frac{1}{2}q^{i!}(q^{2i!}-1)<q^{(i+1)!}=|U_{i+1}|. $$
The contradiction between (\ref{35}), (\ref{36}) completes the proof.
\end{proof}
It follows immediately from Lemma \ref{noHG} that (\ref{exth}) is nonsplit, and Proposition \ref{main2} is proved.

\begin{proposition}\label{main3}
If $\theta|_{\bc}={\rm tr}$, then $\op{Ext}_{\bk\bg}^{1}({\rm St},\m(\theta))\ne0$
\end{proposition}

We will construct a nonsplit sequence
 \begin{equation}\label{extl}
0\rightarrow\m(\theta)\rightarrow L \rightarrow \text{St}  \rightarrow0
 \end{equation}
of $\bk\bg$-modules and prove  Proposition \ref{main3}. Noting that  ${\rm St}=\Bbbk {\bf U} \eta$, we denote  ${\rm St}_i= \Bbbk {U_i} \eta$.  Let $L_i={\rm St}_i \oplus M_i(\theta)$ which is a $\bk G_i$-module. For each $i$, we fix an element  $b_i\in\f_{q^{(i+1)!}}-\f_{q^{2i!}}$ and set $$\zeta_i=(1-s) \underline{U_i} \underline{T_i} \varepsilon(b_i)\sn\in M_{i+1}(\theta). $$ By the same arguments  in the proof of  Lemma \ref{ne02}, we have the following lemma.

\begin{lemma}\label{ne03}
One has that $\zeta_i\ne 0$. Moreover,
$$s\zeta_i=-\zeta_i,  \quad s\varepsilon(x) \zeta_i=(\varepsilon(-x^{-1})-1)\zeta_i$$
where $x\in U_i$ and  $t\zeta_i=\zeta_i$ for any $t\in T_i$.
\end{lemma}

\begin{proof}  Firstly, we have
$$\zeta_i= (1-s)\sum_{{t\in T_i}\atop{a\in\f_{q^{i!}}}}\theta(t)^{-1}\varepsilon(b_it^2+a)\sn.$$
 Using $ s\varepsilon(a)s=
\varepsilon(-a^{-1})sh(a)\varepsilon(-a^{-1})$, it is not difficult to get
\begin{equation}\label{59}
\zeta_i=\sum_{{t\in T_i}\atop{a\in\f_{q^{i!}}}}\theta(t)^{-1}\varepsilon(b_it^2+a)\sn- \sum_{{t\in T_i}\atop{a\in\f_{q^{i!}}}} \theta(b_it+at^{-1})\varepsilon(-(b_it^2+a)^{-1})\sn.
\end{equation}
Noting that $b_i\in\f_{q^{(i+1)!}}-\f_{q^{2i!}}$, the elements in
$$ \{ \varepsilon(b_it^2+a), \varepsilon(-(b_it^2+a)^{-1}) \mid  {t\in T_i},  {a\in\f_{q^{i!}}} \}$$
are all pairwise different. In particular, we have $\zeta_i\ne 0$.

In the following, we show $ s\varepsilon(x) \zeta_i=(\varepsilon(-x^{-1})-1)\zeta_i$.  The discussion  of other parts are obvious.
Using  $ s\varepsilon(a)s=
\varepsilon(-a^{-1})sh(a)\varepsilon(-a^{-1})$, we get
  $$ s\varepsilon(x) \zeta_i = s \underline{U_i} \underline{T_i} \varepsilon(b_i)\sn- \varepsilon(-x^{-1})  s \underline{U_i} \underline{T_i} \varepsilon(b_i)\sn $$
by some easy calculations.  The right hand side of the above equation is just  $(\varepsilon(-x^{-1})-1)\zeta_i$. The lemma is proved.
\end{proof}

We define the map $l_i: L_i \rightarrow L_{i+1}$ by
$$l_i(\sum_{x} a_x \varepsilon(x) \eta +m)= \sum_{x} a_x \varepsilon(x) (\eta+\zeta_i) +m,$$
where $m\in M_i(\theta)$. Since $\bk U_i\zeta_i$ is a free $\bk U_i$-module by (\ref{59}), $l_i$ is an injective homomorphism in $\op{Hom}_{G_i}(L_i,L_{i+1})$ by Lemma \ref{ne03}.
For any $i<j$, we denote by $l_{i,j}:~L_i\rightarrow L_j$ the composition $l_{j-1}\circ\cdots\circ l_{i+1}\circ l_i$ which is clearly in $\op{Hom}_{G_i}(L_i,L_j)$. Then $\{L_i,l_{i,j}\}$ is a direct system of vector spaces. Let $L$ be the direct limit of $\{L_i,l_{i,j}\}$, which is a $\Bbbk {\bf G}$-module. Thus we get the short exact sequence (\ref{extl}) by the same arguments as in the proof of Lemma \ref{22}.

\begin{lemma}\label{noLG}
There is no $v\in L-\m(\theta)$ such that $v\in L^{\bt}$.
\end{lemma}

\begin{proof} Suppose that there is $v\in L-\m(\theta)$ such that $v\in L^{\bt}$. Assume that $v$ is represented by $v_i\in H_i$. It is easy to see that $v_i$ has the form $v_i=a\eta +m$, where $a\ne 0$ and $m\in M_i(\theta)$. Note that $l_i(v_i)=a(\eta+\zeta_i)+m\in L_{i+1}^{T_{i+1}}$.   The  cardinality of  the orbit of $T_{i+1}$ acting on $U_{i+1}\backslash\{e\}$ by conjugation is $\displaystyle \frac{q^{(i+1)!}-1}{2}$.
However using  the expression of $\zeta_i$ in (\ref{59}), we see that the cardinality of the following set
\begin{equation}\label{37}
U_i \cup \bigcup_{{t\in T_i}\atop{a\in\f_{q^{i!}}}}\varepsilon(b_it^2+a) \cup  \bigcup_{{t\in T_i}\atop{a\in\f_{q^{i!}}}}\varepsilon(-(b_it^2+a)^{-1}) \ne U_{i+1}.
\end{equation}
is $q^{2i!}$, which is less than $\displaystyle \frac{q^{(i+1)!}-1}{2}$.  This is a contraction to $l_i(v_i)\in L_{i+1}^{T_{i+1}}$ and the lemma is proved.
\end{proof}

It follows immediately from Lemma \ref{noLG} that (\ref{extl}) is nonsplit, and Proposition \ref{main3} is proved.  In conclusion, we have the following theorem of the extensions of simple modules for  ${\bf G}= SL_2(\bar{\mathbb{F}}_q)$.

\begin{theorem} Let ${\bf G}= SL_2(\bar{\mathbb{F}}_q)$ and assume that $ p \neq \op{char}\Bbbk\geq 3$ or $\op{char}\bk=0$.
One has that  $\op{Ext}_{\bk\bg}^{1}(M, \op{tr})=0$ for any $M\in \op{Irr}({\bf G, \bf T})$. If $\theta|_{\bc}\ne {\rm tr}$, then
$$\op{Ext}_{\bk\bg}^{1}(\mathbb{M}(\theta), \op{St})=\op{Ext}_{\bk\bg}^{1}(\op{tr}, \mathbb{M}(\theta))=\op{Ext}_{\bk\bg}^{1}(\op{St},\mathbb{M}(\theta))=0.$$
If $\lambda|_{\bf C}\ne\mu|_{\bf C}$, we have $\op{Ext}_{\bk\bg}^1(\m(\lambda),\m(\mu))=0$.  For any other cases, the extension $\op{Ext}_{\bk\bg}^{1}(M, N)$ is nonzero, where $M, N\in \op{Irr}({\bf G, \bf T})$.
\end{theorem}

\bigskip

\noindent{\bf Statements and Declarations}  The authors  declare that they have no conflict of interests with others.

\medskip

\noindent{\bf Data Availability}  Data sharing not applicable to this article as no datasets were generated or analysed during the current study.

\medskip

\noindent{\bf Acknowledgements} The authors are grateful to  Nanhua Xi for his inspiring encouragements and helpful comments.   This work is sponsored by  NSFC-12101405.

\bigskip

\bibliographystyle{amsplain}

\begin{thebibliography}{10}

\bibitem {BT}
Borel A, Tits J. \textit{Homomorphismes "abstraits" de groupes algebriques simples}, Ann. of Math. (2) 97 (1973), 499-571.


\bibitem {Car}
R. W. Carter, \textit{Finite groups of Lie type: conjugacy classes and complex characters}, Pure Appl. Math. John Wiley and Sons, New York, 1985.



\bibitem {C1}
Xiaoyu Chen, \textit{On the principal representations of reductive groups in defining characteristic}, J. Algebra 620 (2023), 669–689.


\bibitem {C2}
Xiaoyu Chen, \textit{Irreducible modules of reductive groups with Borel-stable line},  arXiv:2011.04115.

\bibitem {C3}
Xiaoyu Chen,  \textit{The principal representation category of infinite reductive groups is not a highest weight category},  Comm. Algebra 51 (2023), no. 1, 157-160.

\bibitem {CD1}
Xiaoyu Chen, Junbin Dong, \textit{Abstract-induced modules for reductive algebraic groups with Frobenius maps},  Int. Math. Res. Not. IMRN 2022, no. 5, 3308-3348.


\bibitem {CPS}
E. Cline,  B. Parshall,  L. Scott,  \textit{Finite-dimensional algebras and highest weight categories}, J. Reine Angew. Math. 391 (1988), 85-99.


\bibitem {D1}
Junbin Dong, \textit{The principal representations of reductive  algebraic groups with Frobenius maps}, J. Algebra 591 (2022), 342-359.


\bibitem {D2} Junbin Dong, \textit{Complex representations of reductive algebraic groups with Frobenius maps in the category $\mathscr{X}$}, arXiv:2309.00341.


\bibitem {D3}  Junbin Dong, \textit{Certain complex representations of $SL_2(\bar{\mathbb{F}}_q)$}, J. Algebra 612 (2022), 504–525.


\bibitem {Xi}
Nanhua Xi,  \textit{Some infinite dimensional representations of reductive groups with Frobenius maps}, Sci. China Math. 57(2014), 1109-1120.

\end{thebibliography}

\end{document}